# ADAPTIVE CONFIDENCE BALLS[1]

By T. Tony Cai and Mark G. Low

*University of Pennsylvania*

Adaptive confidence balls are constructed for individual resolution levels as well as the entire mean vector in a multiresolution framework. Finite sample lower bounds are given for the minimum expected squared radius for confidence balls with a prespecified confidence level. The confidence balls are centered on adaptive estimators based on special local block thresholding rules. The radius is derived from an analysis of the loss of this adaptive estimator. In addition adaptive honest confidence balls are constructed which have guaranteed coverage probability over all of $\mathbb{R}^N$ and expected squared radius adapting over a maximum range of Besov bodies.

**1. Introduction.** A central goal in nonparametric function estimation, and one which has been the focus of much attention in the statistics literature, is the construction of adaptive estimators. Informally, an adaptive procedure automatically adjusts to the smoothness properties of the underlying function. A common way to evaluate such a procedure is to compute its maximum risk over a collection of parameter spaces and to compare these values to the minimax risk over each of them.

It should be stressed that such adaptive estimators do not provide a data-dependent estimate of the loss, nor do they immediately yield easily constructed adaptive confidence sets. Such confidence sets should have size which adapts to the smoothness of the underlying function while maintaining a prespecified coverage probability over a given function space. Moreover, it is clearly desirable to center such confidence sets on estimators which possess other strong optimality properties. In the present paper, a confidence ball is constructed centered on a special block thresholding rule which has particularly good spatial adaptivity. The radius is built upon good estimates of loss.

Received May 2004; revised April 2005.

[1]Supported in part by NSF Grant DMS-03-06576.

We focus on a sequence of statistical models commonly used in the adaptive estimation literature, namely, a multivariate normal model with mean vector corresponding to wavelet coefficients. More specifically, consider the models

$$(1) \qquad y_{j,k} = \theta_{j,k} + \frac{1}{\sqrt{n}} z_{j,k}, \qquad j = 0, 1, \ldots, J-1, k = 1, \ldots, 2^j,$$

where $z_{j,k} \overset{\text{i.i.d.}}{\sim} N(0,1)$ and where it is assumed that $N$ is a function of $n$, $2^J - 1 = N$ and that the mean vector $\theta$ lies in a parameter space $\Theta$. In the present work, confidence balls are constructed over collections of Besov bodies

$$(2) \qquad B_{p,q}^{\beta}(M) = \left\{ \theta : \left( \sum_{j=0}^{J-1} \left( 2^{js} \left( \sum_{k=1}^{2^j} |\theta_{j,k}|^p \right)^{1/p} \right)^q \right)^{1/q} \le M \right\},$$

where $s = \beta + \frac{1}{2} - \frac{1}{p} > 0$ and $p \ge 2$. In particular, these spaces contain as special cases a number of traditional smoothness classes such as Sobolev and Hölder spaces. Although not needed for the development given in this paper, it may be helpful to think of the $\theta_{j,k}$ as wavelet coefficients of a regression function $f$. A confidence ball for the vector $\theta$ then yields a corresponding confidence ball for the regression function $f$. See, for example, [8], where such an approach is taken. Based on the model (1), we introduce new estimates of the loss of block thresholding estimators and use these estimates to construct confidence balls.

In the context of confidence balls, adaptation over a general collection of parameter spaces $\mathcal{C} = \{\Theta_i : i \in I\}$ where $I$ is an index set can be made precise as follows. An adaptive confidence ball guarantees a given coverage probability over the union of these spaces while simultaneously minimizing the maximum expected squared radius over each of the parameter spaces. Write $\mathcal{B}_{\alpha, \Theta}$ for the collection of all confidence balls which have coverage probability of at least $1 - \alpha$ over $\Theta$. Write $r^2(CB, \Theta)$ for the maximum expected squared radius of a confidence ball $CB$ over $\Theta$ and $r_\alpha^2(\Theta)$ for the minimax expected squared radius over confidence balls in $\mathcal{B}_{\alpha, \Theta}$. Then $r_\alpha^2(\Theta)$ is the smallest maximum expected squared radius of confidence balls with guaranteed coverage over $\Theta$. Adaptation over the collection $\mathcal{C}$ can then be defined as follows. Let $\Theta_I = \bigcup_{i \in I} \Theta_i$. A confidence ball $CB \in \mathcal{B}_{\alpha, \Theta_I}$ is called adaptive over $\mathcal{C}$ if for all $i \in I$, $r^2(CB, \Theta_i) \le C_i r_\alpha^2(\Theta_i)$ where $C_i$ are constants not depending on $n$, and we say that adaptation is possible over $\mathcal{C}$ if such a procedure exists.

In a multivariate normal setup as given in the model (1) with $N = n$, Li [11] constructs adaptive confidence balls for the mean vector which have a given coverage over all of $\mathbb{R}^N$. It was shown that under this constraint the squared radius of the ball must, with high probability, be bounded from below by $cn^{-1/4}$ for all choices of the unknown mean vector. Moreover a confidence ball was constructed centered on a shrinkage estimator which attains this lower bound at least for some subsets of $\mathbb{R}^N$.

Hoffmann and Lepski [9] introduce the concept of a random normalizing factor into the study of nonparametric function estimation and used this idea to construct asymptotic confidence balls which adapt over a collection of finitely many parameter spaces. In particular, their results can be used to yield asymptotic confidence balls which adapt over a finite number of Sobolev bodies. Baraud [1] is a further development of both Li [11] and Hoffman and Lepski [9] concentrating on confidence balls which perform well over a finite family of linear subspaces. An honest confidence ball over $\mathbb{R}^N$ was constructed such that the radius adapts with high probability to a given collection of subspaces.

Juditsky and Lambert-Lacroix [10] develop adaptive $L_2$ confidence balls for a function $f$ in a nonparametric regression setup with equally spaced design. The paper used unbiased estimates of risk to construct minimax rate adaptive procedures over Besov spaces. It focused on the asymptotic performance and detailed finite sample results were not given. Robins and van der Vaart [12] use sample splitting to divide the construction of the center and radius of a confidence ball into independent problems and show how to use estimates of quadratic functionals to construct adaptive confidence balls.

In the present paper the focus is on finite sample properties of adaptive confidence balls centered on a special local block thresholding estimator known to have strong adaptivity under mean integrated squared error. The radius is derived from an analysis of the loss of this adaptive estimator. The evaluation of the performance of the resulting confidence ball relies on a detailed understanding of the interplay between these two estimates. Three cases of interest are considered in detail. We first construct confidence balls for the mean vector at individual resolution levels. Then adaptive confidence balls are constructed for all $N$ coefficients over Besov bodies. Finally we consider honest confidence balls over all of $\mathbb{R}^N$ and expected squared radius adapting over a maximum range of Besov bodies.

The paper is organized as follows. Section 2 is focused on constructing confidence balls for the mean vector associated with a single resolution level $j$ in the Gaussian model (1). These confidence balls can be used in a multiresolution study. Finite sample lower bounds are given for the expected squared radius of confidence balls which have a prescribed minimum coverage level over a given Besov body. Bounds are given for the maximum expected squared radius as well as when the mean vector is equal to zero. Confidence balls which have an expected squared radius within a constant factor of both these lower bounds are constructed. We show that the problem is degenerate over a certain range of Besov bodies beyond which full adaptation is possible. Adaptive confidence balls are constructed centered on a block thresholding estimator. The results and ideas given in this section are used as building blocks in the analysis and construction of adaptive confidence balls for all $N$ coefficients in Sections 3 and 4.

The focus of Section 3 is on the construction and analysis of confidence balls with a specified minimal coverage probability over a given Besov body $B_{p,q}^{\beta}(M)$. It is shown that the possible range of adaptation depends on the relationship between the dimension $N$ and the noise level. Adaptive confidence balls are constructed over a maximal range of Besov bodies. These results are markedly different from the bounds derived for adaptive estimation or adaptive confidence intervals.

In Section 4 confidence balls are constructed which have guaranteed coverage probability over all of $\mathbb{R}^N$. This procedure has a number of strong optimality properties. It adapts over a maximal range of Besov bodies over which honest confidence balls can adapt. Moreover, given that the confidence ball has a prespecified coverage probability over $\mathbb{R}^N$, it has maximum expected squared radius within a constant factor of the smallest maximum expected squared radius for all Besov bodies $B_{p,q}^{\beta}(M)$ with $\beta > 0$ and $M \geq 1$.

Proofs are given in Section 5.

## 2. Adaptive confidence balls for a single resolution level.
As mentioned in the Introduction, the mean $\theta_{j,k}$ in the model (1) can be thought of as the $k$th coefficient at level $j$ in a wavelet expansion of a function $f$. The different levels $j$ allow for a multiresolution analysis where the coefficients with small values of $j$ correspond to coarse features and where the coefficients with large values of $j$ correspond to fine features. In this section we first fix a level $j$ and focus not only on estimating the sequence of means at that level but also on constructing honest confidence balls for this set of coefficients.

Confidence balls are constructed which maintain coverage no matter the values of $\theta_{j,k}$ and have an expected radius adapting to these coefficients over a range of Besov bodies. The analysis given in this section also provides insight (as is shown in Sections 3 and 4) into the problem of estimating all the wavelet coefficients across different levels.

In the following analysis, for a given level $j$, write $\theta_j$ for the sequence of mean values at this given resolution level. That is, $\theta_j = \{\theta_{j,k} : k = 1, \ldots, 2^j\}$. The analysis can then naturally be divided into two parts. We start with lower bounds for the expected squared radius of confidence balls which have a given coverage probability over a given Besov body. Two lower bounds are given. One is for the expected squared radius when all the coefficients are zero. The other is for the maximum expected squared radius. Set $z_\alpha = \Phi^{-1}(1-\alpha)$, where $\Phi$ is the cumulative distribution function of a standard Normal random variable.

THEOREM 1. *Fix $0 < \alpha < \frac{1}{2}$ and let $CB(\delta, r_\alpha) = \{\theta_j : \|\theta_j - \delta\|_2 \leq r_\alpha\}$ be a confidence ball for $\theta_j$ with random radius $r_\alpha$ which has a guaranteed coverage probability over $B_{p,q}^\beta(M)$ of at least $1 - \alpha$. Then for any $0 < \varepsilon < \frac{1}{2}(\frac{1}{2} - \alpha)$*

$$(3) \qquad \sup_{\theta \in B_{p,q}^\beta(M)} E_\theta(r_\alpha^2) \geq \frac{\varepsilon^2}{1 - \alpha - \varepsilon} \min(M^2 2^{-2\beta j}, z_{\alpha + 2\varepsilon}^2 2^j n^{-1}).$$

*Moreover, for any $0 < \varepsilon < \frac{1}{2} - \alpha$,*

$$(4) \qquad E_0(r_\alpha^2) \geq \tfrac{1}{4}(1 - 2\alpha - 2\varepsilon) \min(M^2 2^{-2\beta j}, \log^{1/2}(1 + \varepsilon^2) 2^{j/2} n^{-1}),$$

*where $E_0$ denotes expectation under $\theta = 0$.*

It is useful to note that the maximum value of $\sum_k \theta_{j,k}^2$ at a given level $j$ over the Besov body $B_{p,q}^\beta(M)$ is $M^2 2^{-2\beta j}$. Hence, from (4), if $M^2 2^{-2\beta j} < \log^{1/2}(1 + \varepsilon^2) 2^{j/2} n^{-1}$ the lower bound for the expected squared radius when the mean vector is equal to zero is a constant multiple of $M^2 2^{-2\beta j}$. It follows that if a given coverage probability is guaranteed over $B_{p,q}^\beta(M)$ then the maximum expected squared radius over any other Besov body must also be of this same order. It should be stressed that this is really a degenerate case since the trivial ball centered at zero with squared radius equal to $M^2 2^{-2\beta j}$ is within a constant factor of the lower bounds given in (3) and (4) and has coverage probability equal to one. Thus we shall focus only on the construction of confidence balls which have a given coverage probability at least over Besov bodies where $M^2 2^{-2\beta j} > \log^{1/2}(1 + \varepsilon^2) 2^{j/2} n^{-1}$. In particular, we only need to consider resolution levels $j = j_n$ satisfying $2^j \leq n^2$ since resolution levels with $2^j > n^2$ satisfy $M^2 2^{-2\beta j} \leq \log^{1/2}(1 + \varepsilon^2) 2^{j/2} n^{-1}$ at least for large $n$. Moreover, since little is to be gained for levels where $2^j \leq \log n$, by using confidence balls with random radius in such cases we shall just use the usual $100(1 - \alpha)\%$ confidence ball centered on the observations $y_{j,k}$. Thus in the following construction attention is focused on cases where $\log n \leq 2^j \leq n^2$.

As mentioned in the Introduction, the center of the ball is constructed by local thresholding. Set $L = \log n$ and let $B_i^j = \{(j, k) : (i - 1)L + 1 \leq k \leq iL\}$, $1 \leq i \leq 2^j/L$, denote the set of indices of the coefficients in the $i$th block at level $j$. For a given block $B_i^j$, set

$$(5) \qquad S_{j,i}^2 \equiv \sum_{(j,k) \in B_i^j} y_{j,k}^2, \qquad \xi_{j,i}^2 \equiv \sum_{(j,k) \in B_i^j} \theta_{j,k}^2 \quad \text{and} \quad \chi_{j,i}^2 \equiv \sum_{(j,k) \in B_i^j} z_{j,k}^2.$$

Let $\lambda_* = 6.9368$ be the root of the equation $\lambda - \log \lambda = 5$. This threshold is similar to the one used in [4, 5]. Then the center $\hat{\theta} = (\hat{\theta}_{j,k})$ is defined by

$$(6) \qquad \hat{\theta}_{j,k} = y_{j,k} \cdot I(S_{j,i}^2 \geq \lambda_* L n^{-1}).$$

It follows from [5] that this local block thresholding rule has strong adaptivity under both global and local risk measures. We now show how the loss $\|\hat{\theta}_j - \theta_j\|_2^2$ of this estimator can be estimated and used in the construction of the radius of the confidence ball. Note that $\theta_{j,k}$ equals either 0 or $y_{j,k}$ and hence the loss can be broken into two terms,

$$
\begin{aligned}
(7) \qquad \sum_k (\hat{\theta}_{j,k} - \theta_{j,k})^2 &= \sum_i \xi_{j,i}^2 I(S_{j,i}^2 \leq \lambda_* L n^{-1}) \\
&\quad + \sum_i n^{-1} \chi_{j,i}^2 I(S_{j,i}^2 > \lambda_* L n^{-1}).
\end{aligned}
$$

The first term can be handled by using an estimate of a quadratic functional. The other term can be analyzed using the fact that $\chi_{j,i}^2$ has a central chi-squared distribution.

Let $(x)_+$ denote $\max(0, x)$ and set

$$
\begin{aligned}
(8) \qquad r_\alpha^2 &= \left[ 2 \log^{1/2}\!\left( \frac{2}{\alpha} \right) + 4 \lambda_*^{1/2} z_{\alpha/2} \right] 2^{j/2} n^{-1} \\
&\quad + \left( \sum_i (S_{j,i}^2 - L n^{-1}) I(S_{j,i}^2 \leq \lambda_* L n^{-1}) \right)_+ \\
&\quad + (2\lambda_* + 8\lambda_*^{1/2} - 1) L n^{-1} \operatorname{Card}\{i : S_{j,i}^2 > \lambda_* L n^{-1}\}.
\end{aligned}
$$

The confidence ball is then defined as

$$
(9) \qquad CB_*(\hat{\theta}_j, r_\alpha) = \{\theta_j : \|\theta_j - \hat{\theta}_j\|_2 \leq r_\alpha\}
$$

where, when $2^j \geq \log n$, the center $\hat{\theta}_j$ is given as in (6) and the radius given in (8) and where $\hat{\theta}_j = y_{j,k}$ and $r_\alpha$ is the radius of the usual $100(1-\alpha)\%$ confidence ball when $2^j < \log n$.

THEOREM 2. *Let the confidence ball $CB_*(\hat{\theta}, r_\alpha)$ be given as in* (9) *and suppose that the resolution level $j$ satisfies $2^j \leq n^2$. Then*

$$
(10) \qquad \inf_{\theta \in \mathbb{R}^N} P(\theta_j \in CB_*(\hat{\theta}, r_\alpha)) \geq 1 - \alpha - 2(\log n)^{-1},
$$

*and for a constant $C_\beta$ depending only on $\beta$,*

$$
\begin{aligned}
(11) \qquad \sup_{\theta \in B_{p,q}^\beta(M)} E(r_\alpha^2) &\leq \left[ 2 \log^{1/2}\!\left( \frac{2}{\alpha} \right) + 4 \lambda_*^{1/2} z_{\alpha/2} + 4 \right] 2^{j/2} n^{-1} \\
&\quad + C_\beta \min(2^j n^{-1}, M^2 2^{-2\beta j}).
\end{aligned}
$$

Note that the confidence ball constructed above attains the minimax lower bound given in (3) simultaneously over all Besov bodies $B_{p,q}^\beta(M)$ with $M^2 2^{-2\beta j} > \log^{1/2}(1 + \varepsilon^2) 2^{j/2} n^{-1}$. This is true even though the confidence ball has a given level of coverage for all $\theta$ in $\mathbb{R}^N$.

**3. Adaptive confidence balls over Besov bodies.** The confidence balls constructed in Section 2 focused on a given resolution level. In this section this construction is extended to the more complicated case of estimating all $N$ coefficients of $\theta$. Specifically, we consider adaptation over a collection of Besov bodies $B_{p,q}^{\beta}(M)$ with $p \geq 2$. It should be stressed that the theory developed in this section for adaptive confidence balls is quite different from that of adaptive estimation theory where adaptation under global losses is possible over all Besov bodies. In particular, adaptation for confidence balls is only possible over a much smaller range of Besov bodies.

In Section 3.1 a lower bound is given on both the maximum and the minimum expected squared radius for any confidence ball with a particular coverage probability over a Besov body. As in Section 2, these lower bounds provide a fundamental limit to the range of Besov bodies where adaptation is possible. Adaptive confidence balls are described in Section 3.2. They build on the construction given in Section 2. The center uses the special local block thresholding rule used in Section 2 up to a particular level and then estimates the remaining coordinates by zero. The radius is chosen based on an estimate of the loss of this block thresholding estimate. The analysis of the resulting confidence ball relies on a detailed understanding of the interplay between these two estimates.

3.1. *Lower bounds.* Theorem 1 provides lower bounds for the expected squared radius of a confidence ball for the mean vector at a given resolution level with a given coverage over $B_{p,q}^{\beta}(M)$. In this section lower bounds are given for the expected squared radius for the whole mean vector for any confidence ball which has a given coverage probability over $B_{p,q}^{\beta}(M)$. There are two lower bounds, one for the maximum expected squared radius and one for the minimum expected squared radius. We shall show that these two lower bounds determine the range over which adaptation is possible.

THEOREM 3. *Fix $0 < \alpha < \frac{1}{2}$ and let $CB(\delta, r_\alpha) = \{\theta : \|\theta - \delta\|_2 \leq r_\alpha\}$ be a $1 - \alpha$ level confidence ball for $\theta \in B_{p,q}^{\beta}(M)$ with random radius $r_\alpha$. Then*

$$
\begin{aligned}
(12) \quad &\sup_{\theta \in B_{p,q}^{\beta}(M)} E_\theta(r_\alpha^2) \\
&\geq \frac{\varepsilon^2}{1 - \alpha - \varepsilon} z_{\alpha+2\varepsilon}^2 \min(Nn^{-1}, z_{\alpha+2\varepsilon}^{-2q/(1+2\beta)} M^{2/(1+2\beta)} n^{-2\beta/(1+2\beta)}).
\end{aligned}
$$

*For any $0 < \varepsilon < \frac{1}{2} - \alpha$, set $\gamma = \log(1 + \varepsilon^2)$. For $0 < M' < M$ set*

$$
\begin{aligned}
(13) \quad b_\varepsilon = \min(&2^{-1/(2(1+4\beta))-1} \gamma^{\beta/(1+4\beta)} (M - M')^{1/(1+4\beta)} n^{-2\beta/(1+4\beta)}, \\
&\tfrac{1}{2} \gamma^{1/4} N^{1/4} n^{-1/2}).
\end{aligned}
$$

*Then for all $\theta \in B_{p,q}^{\beta}(M')$,*

$$
(14) \qquad P_\theta(r_\alpha > b_\varepsilon) \geq 1 - 2\alpha - 2\varepsilon
$$

*and consequently*

$$\inf_{\theta \in B_{p,q}^\beta(M')} E_\theta(r_\alpha^2) \geq (1 - 2\alpha - 2\varepsilon) b_\varepsilon^2. \tag{15}$$

In fact, as is shown in the next section, both bounds are rate sharp in the sense that there are confidence balls with a given coverage probability over $B_{p,q}^\beta(M)$ which have expected squared radius within a constant factor of the lower bounds given in (12) and (15). There are two cases of interest, namely, when $N \geq n^2$ and $N < n^2$. First suppose that $N \geq n^2$ and fix a Besov body $B_{p,q}^\beta(M)$ over which it is assumed that the confidence ball has a given coverage probability. Then by (15) the minimum expected squared radius is at least of order $n^{-4\beta/(1+4\beta)}$. Since from (12) the minimax expected squared radius for confidence balls over $B_{p,q}^\beta(M)$ is of order $n^{-2\tau/(1+2\tau)}$, the confidence ball $CB(\delta, r)$ must have expected squared radius larger than the minimax expected squared radius over any Besov body $B_{p',q'}^\tau(M)$ whenever $\tau > 2\beta$ and $p' \geq 2$. Hence in this case it is impossible to adapt over any Besov body with smoothness index $\tau > 2\beta$. Consequently in this case there is a maximum range of Besov bodies over which full adaptation is possible.

Now suppose that $N < n^2$ and that $N \asymp n^\rho$ where $0 < \rho < 2$. In this case the possible range of adaptation depends on the value of $\rho$. Let $CB(\delta, r)$ be a confidence ball with guaranteed coverage probability over $B_{p,q}^\beta(M)$. First suppose that $\beta \geq \frac{1}{2\rho} - \frac{1}{4}$. Then as above it is easy to check that the minimum expected squared radius is at least of order $n^{-4\beta/(1+4\beta)}$ and that it is impossible to adapt over Besov bodies with $\tau > 2\beta$. On the other hand, suppose that $\beta < \frac{1}{2\rho} - \frac{1}{4}$. Then by (15), the minimum expected squared radius is at least of order $n^{\rho/2-1}$, which is the minimax rate of convergence for the squared radius over a Besov body with $\beta = \frac{1}{\rho} - \frac{1}{2}$. Hence in this case it is impossible to adapt over any Besov body with smoothness index $\tau > \frac{1}{\rho} - \frac{1}{2}$.

In summary, for a confidence ball with a prespecified coverage probability over a Besov body $B_{p,q}^\beta(M)$ the maximum range of Besov bodies $B_{p,q}^\tau(M)$ over which full adaptation is possible is given in Table 1.

TABLE 1

| $N$, $n$ and $\beta$ | Maximum range of adaptation |
|---|---|
| $N \geq n^2$, all $\beta > 0$ | $\beta \leq \tau \leq 2\beta$ |
| $N = n^\rho$ for $0 < \rho < 2$, $\beta \geq \frac{1}{2\rho} - \frac{1}{4}$ | $\beta \leq \tau \leq 2\beta$ |
| $N = n^\rho$ for $0 < \rho < 2$, $0 < \beta \leq \frac{1}{2\rho} - \frac{1}{4}$ | $\beta \leq \tau \leq \frac{1}{\rho} - \frac{1}{2}$ |

3.2. *Construction of adaptive confidence balls.* In this section the focus is on confidence balls which have a given minimal coverage over a particular Besov body. Subject to this constraint, confidence balls are constructed which have expected squared radius adapting across a range of Besov bodies. The resulting balls are shown to be adaptive over the maximal range of Besov bodies given in Table 1 for the first two cases summarized in the table. The third case is covered in Section 4.

The ball is centered on a local thresholding rule and the squared radius is based on an analysis of the loss of this thresholding rule. More specifically, for the center $\hat{\theta}$, let $J_1$ be the largest integer satisfying

$$(16) \qquad 2^{J_1} \leq \min(N, M^{2/(1+2\beta)} n^{1/(1+2\beta)}).$$

For all $j \geq J_1$, set $\hat{\theta}_{j,k} = 0$ and for $j \leq J_1 - 1$ let $\hat{\theta}_{j,k}$ be the local thresholding estimator given in (6). The radius is found by analyzing the loss

$$(17) \qquad \sum_{j=0}^{J-1} \sum_{k=1}^{2^j} (\hat{\theta}_{j,k} - \theta_{j,k})^2 = \sum_{j=0}^{J_1-1} \sum_{k=1}^{2^j} (\hat{\theta}_{j,k} - \theta_{j,k})^2 + \sum_{j=J_1}^{J-1} \sum_{k=1}^{2^j} \theta_{j,k}^2.$$

The first of these terms is handled similarly to that used in (7) and (8). The second component in the loss $\sum_{j=J_1}^{J-1} \sum_{k=1}^{2^j} \theta_{j,k}^2$ is a quadratic functional. It can be estimated well by using an unbiased estimate of $\sum_{j=J_1}^{J_2-1} \sum_{k=1}^{2^j} \theta_{j,k}^2$ where $J_2$ is the largest integer satisfying $2^{J_2} \leq \min(N, M^{4/(1+4\beta)} n^{2/(1+4\beta)})$ and then bounding the tail $\sum_{j=J_2}^{J-1} \sum_{k=1}^{2^j} \theta_{j,k}^2$ from above.

More specifically, set the squared radius

$$
\begin{aligned}
(18) \qquad r_\alpha^2 = {} & c_\alpha M^{2/(1+4\beta)} n^{-4\beta/(1+4\beta)} \\
& + \sum_{j=0}^{J_1-1} \left( \sum_i (S_{j,i}^2 - Ln^{-1}) I(S_{j,i}^2 \leq \lambda_* Ln^{-1}) \right)_+ \\
& + (2\lambda_* + 8\lambda_*^{1/2} - 1) Ln^{-1} \sum_{j=0}^{J_1-1} \sum_i I(S_{j,i}^2 > \lambda_* Ln^{-1}) \\
& + \sum_{j=J_1}^{J_2-1} \sum_{k=1}^{2^j} (y_{j,k}^2 - n^{-1}),
\end{aligned}
$$

where

$$
\begin{aligned}
c_\alpha = {} & 2^{2\beta} (1 - 2^{-2\beta})^{-1} + 2 \log^{1/2}\left(\frac{4}{\alpha}\right) \\
& + \Bigg\{ \left[ 2 \log^{1/2}\left(\frac{4}{\alpha}\right) + z_{\alpha/4} \cdot 2^{5/2} \lambda_*^{1/2} (1 - 2^{-2\beta})^{1/(2+4\beta)} \right. \\
& \qquad\qquad \left. + z_{\alpha/4} \cdot 2^{\beta+1} (1 - 2^{-2\beta})^{-1/2} \right]
\end{aligned}
$$

$$\times M^{1/(1+2\beta)-2/(1+4\beta)}n^{1/(2+4\beta)-1/(1+4\beta)}\Big\}.$$

Note that the last term in $c_\alpha$ tends to 0 as $n \to \infty$ or $M \to \infty$.

The following theorem shows that the confidence ball $CB^*$ defined by

$$(19) \qquad CB^* = \{\theta : \|\theta - \hat{\theta}\|_2 \le r_\alpha\}$$

has adaptive radius and desired coverage probability.

**THEOREM 4.** *Fix $0 < \alpha < \frac{1}{2}$ and let the confidence ball $CB^*$ be given as in (19). Then, for any $\tau \ge \beta$,*

$$
\begin{aligned}
(20) \qquad & \inf_{\theta \in B_{p,q}^\tau(M)} P(\theta \in CB^*) \\
& \ge (1 - \alpha) \\
& \quad - [n^{-1} + 3(1 - 2^{-2\beta})^{-1/(1+2\beta)}]L^{-1}M^{2/(1+2\beta)}n^{-2\beta/(1+2\beta)}.
\end{aligned}
$$

*For $\tau \le 2\beta$,*

$$(21) \qquad \sup_{\theta \in B_{p,q}^\tau(M)} E(r_\alpha^2) \le C_\tau \min(M^{2/(1+2\tau)}n^{-2\tau/(1+2\tau)}, Nn^{-1})$$

*and for $\tau > 2\beta$,*

$$(22) \qquad \sup_{\theta \in B_{p,q}^\tau(M)} E(r_\alpha^2) \le C_\beta \min(M^{2/(1+4\beta)}n^{-4\beta/(1+4\beta)}, Nn^{-1}),$$

*where $C_\tau$ and $C_\beta$ are constants depending only on $\tau$ and $\beta$, respectively.*

Theorem 4 taken together with Theorem 3 shows that the confidence ball $CB^*$ is adaptive over a maximal range of Besov bodies

$$(23) \qquad \mathcal{C} = \{B_{p,q}^\tau(M) : \tau \in [\beta, 2\beta], p \ge 2, q \ge 1\}$$

when either $N \ge n^2$ or $N = n^\rho$, $0 < \rho < 2$ and $\beta \ge \frac{1}{2\rho} - \frac{1}{4}$. In addition, the results also show that the confidence ball $CB^*$ still has guaranteed coverage over $B_{p,q}^\tau(M)$ for $\tau > 2\beta$ although the maximum expected radius is necessarily inflated.

**4. Adaptive confidence balls with coverage over $\mathbb{R}^N$.** In Section 3 it was assumed that the mean vector belongs to a Besov body $B_{p,q}^\beta(M)$ and the confidence ball was constructed to ensure that it had a prespecified coverage probability over that Besov body. Under this constraint there are two situations where the confidence ball has expected squared radius that adapts over the Besov bodies $B_{p,q}^\tau(M)$ with $\tau$ between $\beta$ and $2\beta$, namely, when $N \ge n^2$ or when $N = n^\rho$ with $0 < \rho < 2$ and $\beta \ge \frac{1}{2\rho} - \frac{1}{4}$. In both cases this is the largest range over which adaptation is possible.

We now turn to a construction of "honest" confidence balls which have guaranteed coverage over all of $\mathbb{R}^N$. For the case when $N = n$, such "honest" confidence balls, those with a guaranteed coverage probability over all of $\mathbb{R}^N$, was a topic pioneered in [11]. See also [2] and [3]. Li [11] was the first to show, when $N = n$, that any "honest" confidence ball must have a minimum expected squared radius of order $n^{-1/2}$. In fact, using the lower bounds in Theorem 1 for the level-by-level case, it is easy to see that for any confidence interval with coverage over all of $\mathbb{R}^N$ the random radius must in general satisfy

$$(24) \qquad E_0(r_\alpha^2) \geq \frac{1 - 2\alpha - 2\varepsilon}{4} (\log(1 + \varepsilon^2))^{1/2} \cdot N^{1/2} n^{-1}.$$

Once again, for the case when $N = n$, Li [11] also showed how to construct "honest" confidence balls with maximum expected squared radius of order $n^{-1/2}$ over a parameter space where a linear estimator can be constructed with maximum risk of order $n^{-1/2}$. Such estimators exist when the parameter space only consists of sufficiently smooth functions. In particular, for the Besov bodies $B_{p,q}^\beta(M)$ with $p \geq 2$ Donoho and Johnstone [7] showed that the minimax linear risk is of order $n^{-2\beta/(1+2\beta)}$ and the methodology of Li [11] then leads to "honest" confidence balls with maximum expected squared radius converging at a rate of $n^{-1/2}$ over Besov bodies $B_{p,q}^\beta(M)$ if $\beta \geq \frac{1}{2}$ and $p \geq 2$. However this approach is not adaptive over Besov bodies $B_{p,q}^\beta(M)$ with $\beta \leq \frac{1}{2}$.

In this section "honest" confidence balls are constructed over $\mathbb{R}^N$ which simultaneously adapt over a maximal range of Besov bodies. Attention is focused on the case where $N \leq n^2$ since, from (24), if $N > n^2$, the minimum expected squared radius of such "honest" confidence balls does not even converge to zero.

The confidence ball is built by applying the single level construction given in Section 2 level by level. In particular, the center of the confidence ball is obtained by block thresholding all the observations in blocks of size $L = \log n$. For each index $(j, k)$ in the block, say, $B_i^j$ the estimate of $\theta_{j,k}$ is given by

$$(25) \qquad \hat{\theta}_{j,k} = y_{j,k} \cdot I(S_{j,i}^2 \geq \lambda_* L n^{-1})$$

where $\lambda_* = 6.9368$. The center of the confidence ball $\hat{\theta}$ is then defined by $\hat{\theta} = (\hat{\theta}_{j,k})$. The construction of the radius is once again based on an analysis of the loss $\|\hat{\theta} - \theta\|_2^2$ and applies the same technique as that given in Section

2. Set

$$
\begin{aligned}
(26) \quad r_\alpha^2 =\ & \left[2\log^{1/2}\!\left(\frac{2}{\alpha}\right) + 4\lambda_*^{1/2} z_{\alpha/2}\right] N^{1/2} n^{-1} \\
& + \sum_{j=1}^{J-1}\left(\sum_i (S_{j,i}^2 - Ln^{-1}) I(S_{j,i}^2 \le \lambda_* Ln^{-1})\right)_+ \\
& + (2\lambda_* + 8\lambda_*^{1/2} - 1) Ln^{-1}\,\mathrm{Card}\{i : S_{j,i}^2 > \lambda_* Ln^{-1}\}.
\end{aligned}
$$

With $\hat\theta$ given in (25) and $r_\alpha$ given in (26) the confidence ball is then defined by

$$
(27) \qquad CB_*(\hat\theta, r_\alpha) = \{\theta : \|\theta - \hat\theta\|_2 \le r_\alpha\}.
$$

THEOREM 5. *Let the confidence ball $CB_*(\hat\theta, r_\alpha)$ be given as in* (27). *Then*

$$
(28) \qquad \inf_{\theta \in \mathbb{R}^N} P(\theta \in CB_*(\hat\theta, r_\alpha)) \ge 1 - \alpha - 2(\log n)^{-1}
$$

*and, if $M \ge 1$,*

$$
\begin{aligned}
(29) \quad \sup_{\theta \in B_{p,q}^\tau(M)} E(r_\alpha^2) \le\ & \left[2\log^{1/2}\!\left(\frac{2}{\alpha}\right) + 4\lambda_*^{1/2} z_{\alpha/2} + 4\right] N^{1/2} n^{-1} \\
& + C_\tau \min(Nn^{-1}, M^{2/(1+2\tau)} n^{-2\tau/(1+2\tau)}),
\end{aligned}
$$

*where $C_\tau > 0$ is a constant depending only on $\tau$.*

It is also interesting to understand Theorem 5 from an asymptotic point of view. Fix $0 < \rho < 2$ and let $N = n^\rho$. It then follows from Theorem 5 that the confidence ball constructed above has adaptive squared radius over Besov bodies $B_{p,q}^\tau(M)$ with $\tau \le \frac{1}{\rho} - \frac{1}{2}$ and has maximum expected squared radius of order $n^{-1/2}$ over Besov bodies with $\tau > \frac{1}{\rho} - \frac{1}{2}$. Note that the range depends on $N$. In particular, consider the special case of $N = n$. In this case, note that for $\tau \le \frac{1}{2}$ and $M \ge 1$ it follows that

$$
(30) \qquad \sup_{\theta \in B_{p,q}^\tau(M)} E(r_\alpha^2) \le C_\tau \min(1, M^{2/(1+2\tau)} n^{-2\tau/(1+2\tau)})
$$

and hence, although the confidence ball $CB_*$ depends only on $n$ and the confidence level, it adapts over the collection of all Besov bodies $B_{p,q}^\beta(M)$ with $\beta \le \frac{1}{2}$,

$$
(31) \qquad \mathcal{C} = \{B_{p,q}^\beta(M) : 0 < \beta \le \tfrac{1}{2}, p \ge 2, q \ge 1, M \ge 1\}.
$$

This is the maximal range of Besov bodies over which honest confidence balls can adapt. In addition, it follows from (29) that the confidence ball has maximum expected squared radius within a constant factor of the smallest maximum expected squared radius for all Besov bodies $B_{p,q}^\beta(M)$ with $\beta > 0$ and $M \ge 1$ among all confidence balls which have a prespecified coverage probability over $\mathbb{R}^N$.

**5. Proofs.** In this section proofs of the main theorems are given except for Theorem 2. The proof of Theorem 2 is analogous although slightly easier than that given for Theorem 4.

5.1. *Proof of Theorems 1 and 3.* Theorems 1 and 3 give lower bounds for the squared radius of the confidence balls. A unified proof of these two theorems can be given. We begin with a lemma on the minimax risk over a hypercube.

LEMMA 1. *Suppose $y_i = \theta_i + \sigma z_i$, $z_i \overset{i.i.d.}{\sim} N(0,1)$ and $i = 1, \ldots, m$. Let $a > 0$, and set $C_m(a) = \{\theta \in \mathbb{R}^m : \theta_i = \pm a, i = 1, \ldots, m\}$. Let the loss function be*

$$(32) \qquad L(\hat{\theta}, \theta) = \sum_{i=1}^{m} I(|\hat{\theta}_i - \theta_i| \geq a).$$

*Then the minimax risk over $C_m(a)$ satisfies*

$$(33) \qquad \begin{aligned} \inf_{\hat{\theta}} \sup_{\theta \in C_m(a)} E(L(\hat{\theta}, \theta)) &= \inf_{\hat{\theta}} \sup_{\theta \in C_m(a)} \sum_{i=1}^{m} P(|\hat{\theta}_i - \theta_i| \geq a) \\ &= \Phi\left(-\frac{a}{\sigma}\right) m, \end{aligned}$$

*where $\Phi(\cdot)$ is the cumulative distribution function for the standard Normal distribution.*

PROOF. Let $\pi_i$, $i = 1, \ldots, m$, be independent with $\pi_i(a) = \pi_i(-a) = \frac{1}{2}$. Let $\pi = \prod_{i=1}^{m} \pi_i$ be the product prior on $\theta \in C_m(a)$. The posterior distribution of $\theta$ given $y$ can be easily calculated as $P_{\theta|y}(\theta) = \prod_{i=1}^{m} P_{\theta_i|y_i}(\theta_i)$ where

$$P_{\theta_i|y_i}(\theta_i) = \frac{e^{2ay_i/\sigma^2}}{1 + e^{2ay_i/\sigma^2}} \cdot I(\theta_i = a) + \frac{1}{1 + e^{2ay_i/\sigma^2}} \cdot I(\theta_i = -a).$$

The Bayes estimator $\hat{\theta}^{\pi}$ under the prior $\pi$ and loss $L(\cdot, \cdot)$ given in (32) is then the minimizer of $E_{\theta|y} L(\hat{\theta}, \theta) = \sum_{i=1}^{m} P_{\theta|y}(|\hat{\theta}_i - \theta_i| \geq a)$. A solution is then given by the simple rule $\hat{\theta}_i^{\pi} = a$ if $y_i \geq 0$, $\hat{\theta}_i^{\pi} = -a$ if $y_i < 0$. The risk of the Bayes rule $\hat{\theta}^{\pi}$ equals

$$(34) \qquad \begin{aligned} &\sum_{i=1}^{m} P_{\theta}(|\hat{\theta}_i^{\pi} - \theta_i| \geq a) \\ &= m \cdot \left\{ \frac{1}{2} P(y_i < 0 | \theta_i = a) + \frac{1}{2} P(y_i \geq 0 | \theta_i = -a) \right\} \\ &= \Phi\left(-\frac{a}{\sigma}\right) m. \end{aligned}$$

Since the risk of the Bayes rule $\hat{\theta}^\pi$ is a constant, it equals the minimax risk. $\square$

The proofs of Theorems 1 and 3 are also based on a bound on the $L_1$ distance between a multivariate normal distribution with mean 0 and a mixture of normal distributions with means supported on the union of vertices of a collection of hyperrectangles. Let $C(a, k)$ be the set of $N$-dimensional vectors of which the first $k$ coordinates are equal to $a$ or $-a$ and the remaining coordinates are equal to 0. Then $\mathrm{Card}(C(a, k)) = 2^k$. Let $P_k$ be the mixture of Normal distributions with mean supported over $C(a, k)$,

$$(35) \qquad P_k = \frac{1}{2^k} \sum_{\theta \in C(a,k)} \Phi_{\theta, 1/\sqrt{n}, N},$$

where $\Phi_{\theta, \sigma, N}$ is the Normal distribution $N(\theta, \sigma^2 I_N)$. Denote by $\phi_{\theta, \sigma, N}$ the density of $\Phi_{\theta, \sigma, N}$ and set $P_0 = \Phi_{0, 1/\sqrt{n}, N}$.

LEMMA 2. *Fix $0 < \varepsilon < 1$ and suppose $ka^4 n^2 \leq \log(1 + \varepsilon^2)$. Then*

$$(36) \qquad L_1(P_0, P_k) \leq \varepsilon.$$

*In particular, if $A$ is any event such that $P_0(A) \geq \alpha$, then*

$$(37) \qquad P_k(A) \geq \alpha - \varepsilon,$$

*where $P_k$ is the mixture of Normal distributions given in* (35).

PROOF. The chi-squared distance between the distributions $P_k$ and $P_0 = \Phi_{0, 1/\sqrt{n}, N}$ satisfies $\int \frac{P_k^2}{P_0} \leq e^{ka^4 n^2} \leq 1 + \varepsilon^2$ and consequently the $L_1$ distance between $P_0$ and $P_k$ satisfies

$$L_1(P_0, P_k) = \int |dP_0 - dP_k| \leq \left( \int \frac{P_k^2}{P_0} - 1 \right)^{1/2} \leq \varepsilon.$$

Hence, if $P_0(A) \geq \alpha$, then $P_k(A) \geq P_0(A) - L_1(P_0, P_k) \geq \alpha - \varepsilon$ and the lemma follows. $\square$

PROOF OF THEOREMS 1 AND 3. We first prove the bound (3). Fix a constant $\varepsilon$ satisfying $0 < \varepsilon < \frac{1}{2}(\frac{1}{2} - \alpha)$ and note that $z_{\alpha + 2\varepsilon} > 0$. Take $m = 2^j$, $\sigma = n^{-1/2}$ and $a = \min(z_{\alpha + 2\varepsilon} n^{-1/2}, M 2^{-j(\beta + 1/2)})$ in Lemma 1 and let $C_m(a)$ be defined as in Lemma 1. Then every $N$-dimensional vector with the $j$th level coordinates $\theta_j$ in $C_m(a)$ and other coordinates equal to zero is contained in $B_{p,q}^\beta(M)$. It then follows from Lemma 1 that

$$(38) \quad \inf_{\hat{\theta}} \sup_{\theta \in B_{p,q}^\beta(M)} \sum_{k=1}^m P(|\hat{\theta}_{j,k} - \theta_{j,k}| \geq a) \geq \inf_{\hat{\theta}} \sup_{\theta_j \in C_m(a)} \sum_{k=1}^m P(|\hat{\theta}_{j,k} - \theta_{j,k}| \geq a)$$
$$\geq (\alpha + 2\varepsilon)m.$$

For any $\hat{\theta}$, set $X_\theta = \sum_{k=1}^{m} I(|\hat{\theta}_{j,k} - \theta_{j,k}| \geq a)$. Then $X_\theta \leq m$. Let $\gamma = \frac{\varepsilon}{1-\alpha-\varepsilon}$. Then

$$(\alpha + 2\varepsilon)m \leq \sup_{\theta \in B_{p,q}^\beta(M)} E(X_\theta) \leq \sup_{\theta \in B_{p,q}^\beta(M)} \{\gamma m P(X_\theta < \gamma m) + m P(X_\theta \geq \gamma m)\}.$$

It follows that $\sup_{\theta \in B_{p,q}^\beta(M)} P(X_\theta \geq \gamma m) \geq \alpha + \varepsilon$ and consequently

$$(39) \qquad \sup_{\theta \in B_{p,q}^\beta(M)} P(\|\hat{\theta}_j - \theta_j\|_2^2 \geq \gamma m a^2) \geq \sup_{\theta \in B_{p,q}^\beta(M)} P(X_\theta \geq \gamma m) \geq \alpha + \varepsilon.$$

Suppose $CB(\hat{\theta}, r_\alpha) = \{\theta_j : \|\theta_j - \hat{\theta}_j\|_2 \leq r_\alpha\}$ is a $1 - \alpha$ level confidence ball over $B_{p,q}^\beta(M)$. Then $\inf_{\theta \in B_{p,q}^\beta(M)} P(\|\theta_j - \hat{\theta}\|_2^2 \leq r_\alpha^2) \geq 1 - \alpha$ and hence

$$\sup_{\theta \in B_{p,q}^\beta(M)} P(r_\alpha^2 \geq \gamma m a^2) \geq \sup_{\theta \in B_{p,q}^\beta(M)} P(\gamma m a^2 \leq \|\theta_j - \hat{\theta}_j\|_2^2 \leq r_\alpha^2)$$
$$\geq \alpha + \varepsilon + 1 - \alpha - 1 = \varepsilon.$$

Thus for any $\varepsilon$ satisfying $0 < \varepsilon < \frac{1}{2}(\frac{1}{2} - \alpha)$, $\sup_{\theta \in B_{p,q}^\beta(M)} E(r_\alpha^2) \geq \varepsilon \gamma m a^2$, which completes the proof of (3). The proof of (12) is quite similar. Let $j'$ be the largest integer satisfying $2^{j'} \leq \min(N, (1-2^{-q(\beta+1/2)})^{2/(q(1+2\beta))} z_{\alpha+2\varepsilon}^{-2q/(1+2\beta)} \times M^{2/(1+2\beta)} n^{1/(1+2\beta)})$. Equation (12) in Theorem 3 follows from Lemma 1 by taking $m = 2^{j'}$, $\sigma = n^{-1/2}$ and $a = z_{\alpha+2\varepsilon} n^{-1/2}$.

We now turn to the proof of (4) and (15). For (4) apply Lemma 2 with $k = 2^j$ and $a = \min(M 2^{-j(\beta+1/2)}, \gamma^{1/4} 2^{-j/4} n^{-1/4})$. It is easy to check by using the first term in the minimum that $2^{js} 2^{j/p} a \leq M$. Hence the sequence which is equal to $a$ or $-a$ on the $j$th level and otherwise zero satisfies the Besov constraint (2). Moreover, using the second term in the minimum, it is clear that $k a^4 n^2 \leq \gamma$. For (15) the above remarks hold with $j$ replaced by $J$ and it is clear that the collection $C(a, k)$ of all such sequences is contained in $B_{p,q}^\beta(M)$. It then follows from Lemma 2 that, for $P_k$ defined by (35), $L_1(P_0, P_k) \leq \varepsilon$ and so

$$(40) \qquad P_k(0 \in CB(\delta, r_\alpha)) \geq 1 - \alpha - \varepsilon.$$

Now since for all $\theta \in C(a, k)$, $P(\theta \in CB(\delta, r_\alpha)) \geq 1 - \alpha$ and hence $P(\{C(a, k) \cap CB(\delta, r_\alpha) \neq \varnothing\}) \geq 1 - \alpha$, it follows that

$$(41) \qquad P_k(\{C(a, k) \cap CB(\delta, r_\alpha) \neq \varnothing\}) \geq 1 - \alpha.$$

The Bonferroni inequality applied to equations (40) and (41) then yields

$$(42) \qquad P_k(0 \in CB(\delta, r_\alpha) \cap \{C(a, k) \cap CB(\delta, r_\alpha) \neq \varnothing\}) \geq 1 - 2\alpha - \varepsilon.$$

Once again, since $L_1(P_0, P_k) \leq \varepsilon$ it follows that

$$(43) \qquad P_0(0 \in CB(\delta, r_\alpha) \cap \{C(a, k) \cap CB(\delta, r_\alpha) \neq \varnothing\}) \geq 1 - 2\alpha - 2\varepsilon.$$

Now note that for all $\theta \in C(a, k)$, $\|\theta\|_2 = ak^{1/2} = 2b_\varepsilon$. Hence, if $CB(\delta, r_\alpha)$ contains both 0 and some point $\theta \in C(a, k)$, it follows that the radius $r_\alpha \geq \frac{1}{2}\|\theta\|_2 = b_\varepsilon$ and consequently

$$P_0(r_\alpha > b_\varepsilon) \geq P_0(0 \in CB(\delta, r_\alpha) \cap \{C(a, k) \cap CB(\delta, r_\alpha) \neq \varnothing\}) \geq 1 - 2\alpha - 2\varepsilon.$$

$\square$

5.2. *Proof of Theorem* 4. The proof of Theorem 4 is involved. We first collect in the following lemmas some preparatory results on the tails of chi-squared distributions and Besov bodies. The proofs of these lemmas are straightforward and is thus omitted here. See [6] for detailed proofs.

LEMMA 3. *Let $X_m$ be a random variable having a central chi-squared distribution with $m$ degrees of freedom. If $d > 0$, then*

$$(44) \qquad P(X_m \geq (1 + d)m) \leq \tfrac{1}{2}e^{-(m/2)(d - \log(1+d))}$$

*and consequently $P(X_m \geq (1 + d)m) \leq \tfrac{1}{2}e^{-(1/4)d^2 m + (1/6)d^3 m}$. If $0 < d < 1$, then*

$$(45) \qquad P(X_m \leq (1 - d)m) \leq e^{-(1/4)d^2 m}.$$

LEMMA 4. *Let $y_i = \theta_i + \sigma z_i$, $i = 1, 2, \ldots, L$, $z_i \overset{i.i.d.}{\sim} N(0, 1)$ and let $\lambda_* = 6.9368$ be the constant satisfying $\lambda - \log \lambda = 5$.*

(i) *For $\tau > 0$ let $\lambda_\tau > 1$ denote the constant satisfying $\lambda - \log \lambda = 1 + \frac{4\tau}{1+2\tau}$. If $\sum_{i=1}^{L} \theta_i^2 \leq (\sqrt{\lambda_*} - \sqrt{\lambda_\tau})^2 L\sigma^2$, then*

$$(46) \qquad P\left(\sum_{i=1}^{L} y_i^2 \geq \lambda_* L\sigma^2\right) \leq P\left(\sum_{i=1}^{L} z_i^2 \geq \lambda_\tau L\right) \leq \tfrac{1}{2}e^{-2\tau/(1+2\tau)L}.$$

(ii) *If $\sum_{i=1}^{L} \theta_i^2 \geq 4\lambda_* L\sigma^2$, then*

$$(47) \qquad P\left(\sum_{i=1}^{L} y_i^2 \leq \lambda_* L\sigma^2\right) \leq P\left(\sum_{i=1}^{L} z_i^2 \geq \lambda_* L\right) \leq \tfrac{1}{2}e^{-2L}.$$

LEMMA 5. (i) *For any $\theta \in B_{p,q}^\tau(M)$ and any $0 < m < J - 1$,*

$$(48) \qquad \sum_{j=m}^{J-1} \sum_{k=1}^{2^j} \theta_{j,k}^2 \leq (1 - 2^{-2\tau})^{-1} M^2 2^{-2\tau m}.$$

(ii) *For a constant $a > 0$, set $\mathcal{I} = \{(j, i) : \sum_{(j,k) \in B_i^j} \theta_{j,k}^2 > aLn^{-1}\}$. Then for $p \geq 2$*

$$\sup_{\theta \in B_{p,q}^\tau(M)} \mathrm{Card}(\mathcal{I}) \leq DL^{-1} M^{2/(1+2\tau)} n^{1/(1+2\tau)}, \tag{49}$$

*where $D$ is a constant depending only on $a$ and $\tau$. In particular, $D$ can be taken as $D = 3(1 - 2^{-2\tau})^{-1/(1+2\tau)} a^{-1/(1+2\tau)}$.*

PROOF OF THEOREM 4. The proof is naturally divided into two parts: expected squared radius and the coverage probability. First recall the notation that for a given block $B_i^j$,

$$S_{j,i}^2 \equiv \sum_{(j,k) \in B_i^j} y_{j,k}^2, \qquad \xi_{j,i}^2 \equiv \sum_{(j,k) \in B_i^j} \theta_{j,k}^2 \quad \text{and} \quad \chi_{j,i}^2 \equiv \sum_{(j,k) \in B_i^j} z_{j,k}^2.$$

We begin with the expected squared radius. Let $\tau \geq \beta$ and suppose $\theta \in B_{p,q}^\tau(M)$. From (18) we have

$$
\begin{aligned}
E_\theta(r_\alpha^2) = {}& c_\alpha M^{2/(1+4\beta)} n^{-4\beta/(1+4\beta)} \\
& + \sum_{j=0}^{J_1-1} E_\theta \left( \sum_i (S_{j,i}^2 - Ln^{-1}) I(S_{j,i}^2 \leq \lambda_* Ln^{-1}) \right)_+ \\
& + (2\lambda_* + 8\lambda_*^{1/2} - 1) Ln^{-1} \sum_{j=0}^{J_1-1} \sum_i P_\theta(S_{j,i}^2 > \lambda_* Ln^{-1}) \\
& + \sum_{j=J_1}^{J_2-1} \sum_{k=1}^{2^j} \theta_{j,k}^2 \\
\equiv {}& G_1 + G_2 + G_3 + G_4.
\end{aligned}
\tag{50}
$$

We begin with the term $G_3$. Let $\lambda_\tau$ be defined as in Lemma 4 and set

$$\mathcal{I}_1 = \{(j, i) : j \leq J_1 - 1, \xi_{j,i}^2 > (\sqrt{\lambda_*} - \sqrt{\lambda_\tau})^2 Ln^{-1}\} \tag{51}$$

and

$$\mathcal{I}_2 = \{(j, i) : j \leq J_1 - 1, \xi_{j,i}^2 \leq (\sqrt{\lambda_*} - \sqrt{\lambda_\tau})^2 Ln^{-1}\}. \tag{52}$$

It then follows from Lemmas 4 and 5 that

$$
\begin{aligned}
\sum_{j=0}^{J_1-1} \sum_i P_\theta(S_{j,i}^2 > \lambda_* Ln^{-1}) = {}& \sum_{(j,i) \in \mathcal{I}_1} P(S_{j,i}^2 > \lambda_* Ln^{-1}) \\
& + \sum_{(j,i) \in \mathcal{I}_2} P(S_{j,i}^2 > \lambda_* Ln^{-1})
\end{aligned}
$$

$$\leq \mathrm{Card}(\mathcal{I}_1) + \tfrac{1}{2} L^{-1} 2^{J_1} \cdot n^{-2\tau/(1+2\tau)}$$

$$\leq \min(L^{-1} 2^{J_1}, D L^{-1} M^{2/(1+2\tau)} n^{1/(1+2\tau)})$$

$$+ \tfrac{1}{2} L^{-1} 2^{J_1} n^{-2\tau/(1+2\tau)}$$

for some constant $D$ depending only on $\tau$. Note that $2^{J_1} = \min(N, M^{2/(1+2\beta)} \times n^{1/(1+2\beta)})$ and so

$$
\begin{aligned}
(53) \quad G_3 &= (2\lambda_* + 8\lambda_*^{1/2} - 1) L n^{-1} \sum_{j=0}^{J_1-1} \sum_i P_\theta(S_{j,i}^2 > \lambda_* L n^{-1}) \\
&\leq C \min(N n^{-1}, M^{2/(1+2\beta)} n^{-2\beta/(1+2\beta)}, M^{2/(1+2\tau)} n^{-2\tau/(1+2\tau)}) \\
&\quad + C \min(N n^{-1}, M^{2/(1+2\beta)} n^{-2\beta/(1+2\beta)}) \cdot n^{-2\tau/(1+2\tau)} \\
&\leq C \min(N n^{-1}, M^{2/(1+2\tau)} n^{-2\tau/(1+2\tau)}).
\end{aligned}
$$

The term $G_4$ is easy to bound. When $N \leq M^{2/(1+2\beta)} n^{1/(1+2\beta)}$, $J_1 = J_2$ and hence $G_4 = 0$. When $N > M^{2/(1+2\beta)} n^{1/(1+2\beta)}$, it follows from (48) in Lemma 5 that

$$
\begin{aligned}
(54) \quad G_4 &= \sum_{j=J_1}^{J_2-1} \sum_{k=1}^{2^j} \theta_{j,k}^2 \\
&\leq (1 - 2^{-2\tau})^{-1} M^2 2^{-2\tau J_1} \\
&\leq C \min(N n^{-1}, M^{2/(1+2\tau)} n^{-2\tau/(1+2\tau)}).
\end{aligned}
$$

We now turn to $G_2$. Let $J_\tau$ be the largest integer satisfying $2^{J_\tau} \leq \min(N, M^{2/(1+2\tau)} n^{1/(1+2\tau)})$. Write

$$
\begin{aligned}
G_2 &= \sum_{j=0}^{J_\tau-1} E\left( \sum_i (S_{j,i}^2 - L n^{-1}) I(S_{j,i}^2 \leq \lambda_* L n^{-1}) \right)_+ \\
&\quad + \sum_{j=J_\tau}^{J_1-1} E\left( \sum_i (S_{j,i}^2 - L n^{-1}) I(S_{j,i}^2 \leq \lambda_* L n^{-1}) \right)_+ \\
&\equiv G_{21} + G_{22},
\end{aligned}
$$

where $G_{22} = 0$ when $J_\tau = J_1$. Note that

$$
\begin{aligned}
(55) \quad G_{21} &= \sum_{j=0}^{J_\tau-1} E\left( \sum_i (S_{j,i}^2 - L n^{-1}) I(S_{j,i}^2 \leq \lambda_* L n^{-1}) \right)_+ \\
&\leq \sum_{j=0}^{J_\tau-1} \sum_i (\lambda_* - 1) L n^{-1} \\
&\leq (\lambda_* - 1) L n^{-1} 2^{J_\tau} L^{-1} \\
&\leq (\lambda_* - 1) \min(N n^{-1}, M^{2/(1+2\tau)} n^{-2\tau/(1+2\tau)}).
\end{aligned}
$$

When $N \leq M^{2/(1+2\tau)} n^{1/(1+2\tau)}$, $J_\tau = J_1$ and so $G_{22} = 0$. On the other hand, when $J_\tau < J_1$,

$$
\begin{aligned}
G_{22} &= \sum_{j=J_\tau}^{J_1-1} E\left(\sum_i (S_{j,i}^2 - Ln^{-1}) I(S_{j,i}^2 \leq \lambda_* Ln^{-1})\right)_+ \\
&\leq \sum_{j=J_\tau}^{J_1-1} E\left(\sum_i (S_{j,i}^2 - Ln^{-1})\right)_+ \\
&\leq \sum_{j=J_\tau}^{J_1-1} \left\{ E\left(\sum_i S_{j,i}^2 - 2^j n^{-1}\right)^2 \right\}^{1/2} \\
&\leq \sum_{j=J_\tau}^{J_1-1} \left\{ 4n^{-1} \sum_i \xi_{j,i}^2 + 2n^{-2} 2^j + \left(\sum_i \xi_{j,i}^2\right)^2 \right\}^{1/2} \\
&\leq 2n^{-1/2} \sum_{j=J_\tau}^{J_1-1} \left(\sum_i \xi_{j,i}^2\right)^{1/2} + 2^{1/2} n^{-1} \sum_{j=J_\tau}^{J_1-1} 2^{j/2} + \sum_{j=J_\tau}^{J_1-1} \sum_i \xi_{j,i}^2.
\end{aligned}
$$

Note that $\sum_i \xi_{j,i}^2 = \sum_{k=1}^{2^j} \theta_{j,k}^2 \leq M^2 2^{-2\tau j}$. It then follows that

$$
\begin{aligned}
G_{22} &\leq 2^{\tau+1}(1 - 2^{-\tau})^{-1} M^{1/(1+2\tau)} n^{-(1+4\tau)/(2+4\tau)} \\
&\quad + 4M^{1/(1+2\beta)} n^{-(1+4\beta)/(2+4\beta)} \\
&\quad + 2^{2\tau}(1 - 2^{-2\tau})^{-1} M^{2/(1+2\tau)} n^{-2\tau/(1+2\tau)}
\end{aligned}
$$

and so

$$
(56) \qquad G_{22} \leq C \min(Nn^{-1}, M^{2/(1+2\tau)} n^{-2\tau/(1+2\tau)}).
$$

This together with (50) and (53)–(55) yields

$$
\begin{aligned}
\sup_{\theta \in B_{p,q}^\tau(M)} E_\theta(r_\alpha^2) &\leq \sup_{\theta \in B_{p,q}^\tau(M)} (G_1 + G_{21} + G_{22} + G_3 + G_4) \\
&\leq c_\alpha \min(Nn^{-1}, M^{2/(1+4\beta)} n^{-4\beta/(1+4\beta)}) \\
&\quad + C_\tau \min(Nn^{-1}, M^{2/(1+2\tau)} n^{-2\tau/(1+2\tau)})
\end{aligned}
$$

where $C_\tau$ is a constant depending only on $\tau$. For $0 < \tau < \beta$ similar arguments yield

$$
\sup_{\theta \in B_{p,q}^\tau(M)} E_\theta(r_\alpha^2) \leq C_\tau \min(Nn^{-1}, M^{2/(1+2\tau)} n^{-2\tau/(1+2\tau)}).
$$

We now turn to the coverage probability. Set $C(\theta) = P(\|\hat{\theta} - \theta\|_2^2 > r_\alpha^2)$ and fix $\tau \geq \beta$. We want to bound $\sup_{\theta \in B_{p,q}^\tau(M)} C(\theta)$. Note that

$$\|\hat{\theta} - \theta\|_2^2 = \sum_{j=0}^{J_1-1} \sum_i \xi_{j,i}^2 I(S_{j,i}^2 \leq \lambda_* L n^{-1})$$

$$+ n^{-1} \sum_{j=0}^{J_1-1} \sum_i \chi_{j,i}^2 I(S_{j,i}^2 > \lambda_* L n^{-1}) + \sum_{j=J_1}^{J-1} \sum_{k=1}^{2^j} \theta_{j,k}^2.$$

It follows from (48) in Lemma 5 that

$$(57) \quad \sup_{\theta \in B_{p,q}^\tau(M)} \sum_{j=J_2}^{J-1} \sum_{k=1}^{2^j} \theta_{j,k}^2 \leq (1 - 2^{-2\tau})^{-1} M^2 2^{-2\tau J_2}$$

$$\leq 2^{2\beta}(1 - 2^{-2\beta})^{-1} M^{2/(1+4\beta)} n^{-4\beta/(1+4\beta)}.$$

Set $a_0 = 2^{2\beta}(1 - 2^{-2\beta})^{-1}$, $a_1 = z_{\alpha/4} \cdot 2^{5/2} \lambda_*^{1/2} (1 - 2^{-2\beta})^{1/(2+4\beta)} \times M^{1/(1+2\beta)-2/(1+4\beta)} n^{1/(2+4\beta)-1/(1+4\beta)}$, $a_2 = 2 \log^{1/2}(\frac{4}{\alpha}) M^{1/(1+2\beta)-2/(1+4\beta)} \times n^{1/(2+4\beta)-1/(1+4\beta)}$, $a_3 = z_{\alpha/4} \cdot 2^{\beta+1}(1 - 2^{-2\beta})^{-1/2} M^{1/(1+2\beta)-2/(1+4\beta)} \times n^{1/(2+4\beta)-1/(1+4\beta)}$, $a_4 = 2 \log^{1/2}(\frac{4}{\alpha})$ and $a_5 = 2\lambda_* + 8\lambda_*^{1/2} - 1$. Then $c_\alpha$ in (18) equals $a_0 + a_1 + a_2 + a_3 + a_4$ and the squared radius $r_\alpha^2$ given in (18) can be written as

$$r_\alpha^2 = (a_0 + a_1 + a_2 + a_3 + a_4) M^{2/(1+4\beta)} n^{-4\beta/(1+4\beta)}$$

$$+ \sum_{j=0}^{J_1-1} \left( \sum_i (S_{j,i}^2 - L n^{-1}) I(S_{j,i}^2 \leq \lambda_* L n^{-1}) \right)_+$$

$$+ a_5 L n^{-1} \sum_{j=0}^{J_1-1} \sum_i I(S_{j,i}^2 > \lambda_* L n^{-1}) + \sum_{j=J_1}^{J_2-1} \sum_{k=1}^{2^j} (y_{j,k}^2 - n^{-1}).$$

Set $\mathcal{I}_3 = \{(j,i) : j \leq J_1 - 1, \xi_{j,i}^2 \geq 4\lambda_* L n^{-1}\}$ and $\mathcal{I}_4 = \{(j,i) : j \leq J_1 - 1, \xi_{j,i}^2 < 4\lambda_* L n^{-1}\}$. It then follows that

$$C(\theta) \leq P\Bigg\{ \sum_{(j,i) \in \mathcal{I}_3} [\xi_{j,i}^2 I(S_{j,i}^2 \leq \lambda_* L n^{-1}) + n^{-1} \chi_{j,i}^2 I(S_{j,i}^2 > \lambda_* L n^{-1})]$$

$$> \sum_{(j,i) \in \mathcal{I}_3} [(S_{j,i}^2 - L n^{-1}) I(S_{j,i}^2 \leq \lambda_* L n^{-1})$$

$$+ a_5 L n^{-1} I(S_{j,i}^2 > \lambda_* L n^{-1})] \Bigg\}$$

$$+ P\Bigg\{ \sum_{(j,i)\in\mathcal{I}_4} [\xi_{j,i}^2 I(S_{j,i}^2 \le \lambda_* L n^{-1}) + n^{-1}\chi_{j,i}^2 I(S_{j,i}^2 > \lambda_* L n^{-1})]$$

$$> (a_1 + a_2) M^{2/(1+4\beta)} n^{-4\beta/(1+4\beta)}$$

$$+ \sum_{(j,i)\in\mathcal{I}_4} [(S_{j,i}^2 - L n^{-1}) I(S_{j,i}^2 \le \lambda_* L n^{-1})$$

$$+ a_5 L n^{-1} I(S_{j,i}^2 > \lambda_* L n^{-1})]\Bigg\}$$

$$+ P\Bigg(\sum_{j=J_1}^{J_2-1} \sum_{k=1}^{2^j} \theta_{j,k}^2 > (a_3 + a_4) M^{2/(1+4\beta)} n^{-4\beta/(1+4\beta)}$$

$$+ \sum_{j=J_1}^{J_2-1} \sum_{k=1}^{2^j} (y_{j,k}^2 - n^{-1})\Bigg)$$

$$\equiv T_1 + T_2 + T_3.$$

We shall consider the three terms separately. We first calculate the term $T_1$. Note that

$$T_1 \le P\Bigg\{ \sum_{(j,i)\in\mathcal{I}_3} (S_{j,i}^2 - \xi_{j,i}^2 - L n^{-1}) I(S_{j,i}^2 \le \lambda_* L n^{-1}) < 0 \Bigg\}$$

$$+ P\Bigg\{ \sum_{(j,i)\in\mathcal{I}_3} n^{-1}\chi_{j,i}^2 I(S_{j,i}^2 > \lambda_* L n^{-1}) > \sum_{(j,i)\in\mathcal{I}_3} a_5 L n^{-1} I(S_{j,i}^2 > \lambda_* L n^{-1}) \Bigg\}$$

$$\le \sum_{(j,i)\in\mathcal{I}_3} P(S_{j,i}^2 \le \lambda_* L n^{-1}) + \sum_{(j,i)\in\mathcal{I}_3} P(\chi_{j,i}^2 > a_5 L).$$

It follows from Lemma 4(ii) that $P(S_{j,i}^2 \le \lambda_* L n^{-1}) \le P(\chi_{j,i}^2 > \lambda_* L) \le \frac{1}{2} n^{-2}$ for $(j,i) \in \mathcal{I}_3$. Lemma 5 now yields

$$(58) \quad \begin{aligned} T_1 &\le n^{-2} \cdot \mathrm{Card}(\mathcal{I}_3) \\ &\le 3(1 - 2^{-2\tau})^{-1/(1+2\tau)} (4\lambda_*)^{-1/(1+2\tau)} L^{-1} M^{2/(1+2\tau)} n^{-2\tau/(1+2\tau)} \\ &\le 3(1 - 2^{-2\beta})^{-1/(1+2\beta)} L^{-1} M^{2/(1+2\beta)} n^{-2\beta/(1+2\beta)}. \end{aligned}$$

We now turn to the second term $T_2$. Note that

$$T_2 = P\Bigg\{ \sum_{(j,i)\in\mathcal{I}_4} [(S_{j,i}^2 - \xi_{j,i}^2 - L n^{-1}) + a_5 L n^{-1} I(S_{j,i}^2 > \lambda_* L n^{-1})]$$

$$< -(a_1 + a_2) M^{2/(1+4\beta)} n^{-4\beta/(1+4\beta)}$$

$$+ \sum_{(j,i) \in \mathcal{I}_4} [(S_{j,i}^2 + n^{-1}\chi_{j,i}^2 - \xi_{j,i}^2 - Ln^{-1})I(S_{j,i}^2 > \lambda_* Ln^{-1})] \Bigg\}$$

$$\leq P \Bigg\{ \sum_{(j,i) \in \mathcal{I}_4} (S_{j,i}^2 - \xi_{j,i}^2 - Ln^{-1}) < -(a_1 + a_2)M^{2/(1+4\beta)}n^{-4\beta/(1+4\beta)} \Bigg\}$$

$$+ P \Bigg\{ \sum_{(j,i) \in \mathcal{I}_4} (S_{j,i}^2 + n^{-1}\chi_{j,i}^2 - \xi_{j,i}^2 - Ln^{-1})I(S_{j,i}^2 > \lambda_* Ln^{-1})$$

$$> \sum_{(j,i) \in \mathcal{I}_4} a_5 Ln^{-1} I(S_{j,i}^2 > \lambda_* Ln^{-1}) \Bigg\}$$

$$\equiv T_{21} + T_{22}.$$

For any given block, write

$$S_{j,i}^2 = \sum_{(j,k) \in B_i^j} (\theta_{j,k} + n^{-1/2}z_{j,k})^2$$

$$= \xi_{j,i}^2 + 2n^{-1/2} \sum_{(j,k) \in B_i^j} \theta_{j,k}z_{j,k} + n^{-1}\chi_{j,i}^2$$

$$= \xi_{j,i}^2 + 2n^{-1/2}\xi_{j,i}\tilde{Z}_{j,i} + n^{-1}\chi_{j,i}^2,$$

where $\tilde{Z}_{j,i} = \xi_{j,i}^{-1} \sum_{(j,k) \in B_i^j} \theta_{j,k}z_{j,k}$ is a standard Normal variable. Then

$$T_{21} = P \Bigg\{ \sum_{(j,i) \in \mathcal{I}_4} (S_{j,i}^2 - \xi_{j,i}^2 - Ln^{-1}) < -(a_1 + a_2)M^{2/(1+4\beta)}n^{-4\beta(1+4\beta)} \Bigg\}$$

$$\leq P \Bigg\{ \sum_{(j,i) \in \mathcal{I}_4} (2n^{-1/2}\xi_{j,i}\tilde{Z}_{j,i} + n^{-1}\chi_{j,i}^2 - Ln^{-1})$$

$$< -(a_1 + a_2)M^{2/(1+4\beta)}n^{-4\beta/(1+4\beta)} \Bigg\}$$

$$\leq P \Bigg\{ 2n^{-1/2} \sum_{(j,i) \in \mathcal{I}_4} \xi_{j,i}\tilde{Z}_{j,i} < -a_1 M^{2/(1+4\beta)}n^{-4\beta/(1+4\beta)} \Bigg\}$$

$$+ P \Bigg\{ \sum_{(j,i) \in \mathcal{I}_4} \chi_{j,i}^2 < -a_2 M^{2/(1+4\beta)}n^{1/(1+4\beta)} + \text{Card}(\mathcal{I}_4)L \Bigg\}$$

$$\equiv T_{211} + T_{212}.$$

Note that, for any $0 < j' \leq J_1 - 1$,

$$\sum_{(j,i) \in \mathcal{I}_4} \xi_{j,i}^2 \leq L^{-1} 2^{j'} \cdot 4\lambda_* L n^{-1} + \sum_{j=j'}^{J_1-1} M^2 2^{-j2\tau}$$

$$= 4\lambda_* 2^{j'} n^{-1} + M^2 (1 - 2^{-2\tau})^{-1} 2^{-2\tau j'}.$$

Minimizing the right-hand side yields that $\sum_{(j,i) \in \mathcal{I}_4} \xi_{j,i}^2 \leq 2(4\lambda_*)^{2\tau/(1+2\tau)} (1 - 2^{-2\tau})^{-1/(1+2\tau)} M^{2/(1+2\tau)} n^{-2\tau/(1+2\tau)} \leq 8\lambda_* (1 - 2^{-2\beta})^{-1/(1+2\beta)} M^{2/(1+2\beta)} \times n^{-2\beta/(1+2\beta)}$. Denote by $Z$ a standard Normal random variable. It then follows that $T_{211} = P(Z < -\frac{1}{2} a_1 M^{2/(1+4\beta)} n^{1/(1+4\beta)} n^{-1/2} (\sum_{(j,i) \in \mathcal{I}_4} \xi_{j,i}^2)^{-1/2}) \leq P(Z < z_{\alpha/4}) = \frac{\alpha}{4}$. Now consider the term $T_{212}$. If $\mathrm{Card}(\mathcal{I}_4) L \leq a_2 M^{2/(1+4\beta)} n^{1/(1+4\beta)}$, then $T_{212} = 0$. Now suppose $\mathrm{Card}(\mathcal{I}_4) L > a_2 M^{2/(1+4\beta)} n^{1/(1+4\beta)}$. It follows from (45) in Lemma 3 by taking $m = \mathrm{Card}(\mathcal{I}_4) L \leq 2^{J_1}$ and $d = a_2 M^{2/(1+4\beta)} \times n^{1/(1+4\beta)}/m$ that $T_{212} \leq \exp(-\frac{1}{4} a_2^2 M^{4/(1+4\beta)-2/(1+2\beta)} n^{2/(1+4\beta)-1/(1+2\beta)}) = \frac{\alpha}{4}$ and hence

(59) $$T_{21} = T_{211} + T_{212} \leq \frac{\alpha}{2}.$$

We now consider the term $T_{22}$. Simple algebra yields that

$$T_{22} = P\Bigg( \sum_{(j,i) \in \mathcal{I}_4} (S_{j,i}^2 + n^{-1} \chi_{j,i}^2 - \xi_{j,i}^2 - L n^{-1}) I(S_{j,i}^2 > \lambda_* L n^{-1})$$

$$> \sum_{(j,i) \in \mathcal{I}_4} a_5 L n^{-1} I(S_{j,i}^2 > \lambda_* L n^{-1}) \Bigg)$$

$$\leq \sum_{(j,i) \in \mathcal{I}_4} P(\tilde{Z}_{j,i} > \tfrac{1}{2} \xi_{j,i}^{-1} (a_5 - 2\lambda_* + 1) L n^{-1/2})$$

$$+ \sum_{(j,i) \in \mathcal{I}_4} P(\chi_{j,i}^2 > \lambda_* L).$$

Note that $\xi_{j,i}^2 \leq 4\lambda_* L n^{-1}$ for $(j,i) \in \mathcal{I}_4$. Hence it follows from the bounds on the tail probability of standard Normal and central chi-squared distributions that

(60) $$T_{22} \leq \sum_{(j,i) \in \mathcal{I}_4} P(\tilde{Z}_{j,i} > 2(\log n)^{1/2}) + \sum_{(j,i) \in \mathcal{I}_4} \tfrac{1}{2} n^{-2}$$
$$\leq L^{-1} M^{2/(1+2\beta)} n^{-2\beta/(1+2\beta)} n^{-1}.$$

We now turn to the third term $T_3$. Note that $y_{j,k}^2 = \theta_{j,k}^2 + 2n^{-1/2} \theta_{j,k} z_{j,k} +$

$n^{-1}z_{j,k}^2$ and so

$$
\begin{aligned}
T_3 = P\bigg( & 2n^{-1/2}\sum_{j=J_1}^{J_2-1}\sum_{k=1}^{2^j}\theta_{j,k}z_{j,k} \\
& + n^{-1}\sum_{j=J_1}^{J_2-1}\sum_{k=1}^{2^j}(z_{j,k}^2-1) < -(a_3+a_4)M^{2/(1+4\beta)}n^{-4\beta/(1+4\beta)}\bigg) \\
\leq P\bigg( & 2n^{1/2}\sum_{j=J_1}^{J_2-1}\sum_{k=1}^{2^j}\theta_{j,k}z_{j,k} < -a_3M^{2/(1+4\beta)}n^{1/(1+4\beta)}\bigg) \\
& + P\bigg(\sum_{j=J_1}^{J_2-1}\sum_{k=1}^{2^j}z_{j,k}^2 < (2^{J_2}-2^{J_1}) - a_4M^{2/(1+4\beta)}n^{1/(1+4\beta)}\bigg) \\
\equiv & T_{31} + T_{32}.
\end{aligned}
$$

Set $\gamma^2 = \sum_{j=J_1}^{J_2-1}\sum_{k=1}^{2^j}\theta_{j,k}^2$ and $Z = \gamma^{-1}\sum_{j=J_1}^{J_2-1}\sum_{k=1}^{2^j}\theta_{j,k}z_{j,k}$. Then $Z$ is a standard Normal variable and it follows from (48) in Lemma 5 that $\gamma^2 \leq 2^{2\beta}(1-2^{-2\beta})^{-1}M^{2/(1+2\beta)}n^{-2\beta/(1+2\beta)}$. Hence,

$$
\begin{aligned}
(61) \quad T_{31} &\leq P(Z < -2^{-\beta-1}(1-2^{-2\beta})^{1/2}a_3M^{2/(1+4\beta)-1/(1+2\beta)} \\
&\qquad\qquad \times n^{1/(1+4\beta)-1/(2+4\beta)}) \\
&= P(Z < -z_{\alpha/4}) = \frac{\alpha}{4}.
\end{aligned}
$$

It follows from Lemma 3 with $m = 2^{J_2} - 2^{J_1}$ and $d = a_4 M^{2/(1+4\beta)}n^{1/(1+4\beta)}/m$ that $T_{32} \leq e^{(-1/4)a_4^2} = \frac{\alpha}{4}$. Equation (20) now follows from this together with (58), (59), (60) and (61). $\quad\square$

5.3. *Proof of Theorem 5.* The proof of Theorem 5 is similar to that of Theorem 4. We shall omit some details and only give a brief proof here. Suppose $\theta \in B_{p,q}^\tau(M)$. Set $b_1 = 2\log^{1/2}(\frac{2}{\alpha})$, $b_2 = 4\lambda_*^{1/2}z_{\alpha/2}$ and $b_3 = 2\lambda_* + 8\lambda_*^{1/2} - 1$. Then, from (26) we have

$$
\begin{aligned}
(62) \quad E_\theta(r_\alpha^2) &= (b_1+b_2)N^{1/2}n^{-1} \\
&\quad + \sum_{j=0}^{J-1}E_\theta\bigg(\sum_i(S_{j,i}^2-Ln^{-1})I(S_{j,i}^2 \leq \lambda_*Ln^{-1})\bigg)_+ \\
&\quad + b_3Ln^{-1}E_\theta(\mathrm{Card}\{(j,i):S_{j,i}^2 > \lambda_*Ln^{-1}\}).
\end{aligned}
$$

The last term can be easily bounded using Lemma 5 as

$$
\begin{aligned}
& b_3Ln^{-1}E_\theta(\mathrm{Card}\{(j,i):S_{j,i}^2 > \lambda_*Ln^{-1}\}) \\
& \qquad \leq b_3\cdot\min(Nn^{-1}, D_\tau M^{2/(1+2\tau)}n^{-2\tau/(1+2\tau)}).
\end{aligned}
$$

Set $D = \sum_{j=0}^{J-1} E_\theta(\sum_i (S_{j,i}^2 - Ln^{-1})I(S_{j,i}^2 \le \lambda_* Ln^{-1}))_+$. Using nearly identical arguments given in the derivation of (55) and (56) in the proof of Theorem 4, $D$ is bounded as $D \le 4N^{1/2}n^{-1} + C_\tau M^{2/(1+2\tau)}n^{-2\tau/(1+2\tau)}$ for some constant $C_\tau \ge \lambda_*$. On the other hand, it is easy to see that $D \le \sum_{j=0}^{J-1} \sum_i (\lambda_* Ln^{-1} - Ln^{-1}) = (\lambda_* - 1)Nn^{-1}$ and consequently $\sup_{\theta \in B_{p,q}^\tau(M)} E(r_\alpha^2) \le (b_1 + b_2 + 4)N^{1/2}n^{-1} + C_\tau \min(Nn^{-1}, M^{2/(1+2\tau)}n^{-2\tau/(1+2\tau)})$.

We now turn to the coverage probability. Again, set $C(\theta) = P(\|\hat\theta - \theta\|_2^2 > r_\alpha^2)$. We want to show that $\sup_{\theta \in \mathbb{R}^N} C(\theta) \le \alpha + 4(\log n)^{-1}$. Note that

$$\|\hat\theta - \theta\|_2^2 = \sum_{j=0}^{J-1} \sum_i \xi_{j,i}^2 I(S_{j,i}^2 \le \lambda_* Ln^{-1})$$
$$+ n^{-1} \sum_{j=0}^{J-1} \sum_i \chi_{j,i}^2 I(S_{j,i}^2 > \lambda_* Ln^{-1}).$$

Set $\mathcal{I}_3' = \{(j,i) : \xi_{j,i}^2 \ge 4\lambda_* Ln^{-1}\}$ and $\mathcal{I}_4' = \{(j,i) : \xi_{j,i}^2 < 4\lambda_* Ln^{-1}\}$. It then follows from the definition of the radius $r_\alpha$ given in (26) that

$$C(\theta) \le P\Bigg\{ \sum_{(j,i) \in \mathcal{I}_3'} [\xi_{j,i}^2 I(S_{j,i}^2 \le \lambda_* Ln^{-1}) + n^{-1}\chi_{j,i}^2 I(S_{j,i}^2 > \lambda_* Ln^{-1})]$$
$$> \sum_{(j,i) \in \mathcal{I}_3'} [(S_{j,i}^2 - Ln^{-1})I(S_{j,i}^2 \le \lambda_* Ln^{-1})$$
$$+ b_3 Ln^{-1} I(S_{j,i}^2 > \lambda_* Ln^{-1})] \Bigg\}$$
$$+ P\Bigg\{ \sum_{(j,i) \in \mathcal{I}_4'} [\xi_{j,i}^2 I(S_{j,i}^2 \le \lambda_* Ln^{-1}) + n^{-1}\chi_{j,i}^2 I(S_{j,i}^2 > \lambda_* Ln^{-1})]$$
$$> (b_1 + b_2)N^{1/2}n^{-1}$$
$$+ \sum_{(j,i) \in \mathcal{I}_4'} [(S_{j,i}^2 - Ln^{-1})I(S_{j,i}^2 \le \lambda_* Ln^{-1})$$
$$+ b_3 Ln^{-1} I(S_{j,i}^2 > \lambda_* Ln^{-1})] \Bigg\}$$
$$\equiv T_1 + T_2.$$

We first bound the term $T_1$. Similarly as in the proof of Theorem 4,

$$(63) \quad T_1 \le \sum_{(j,i) \in \mathcal{I}_3'} (P(S_{j,i}^2 \le \lambda_* Ln^{-1}) + P(\chi_{j,i}^2 > b_3 L)) \le n^{-2} \operatorname{Card}(\mathcal{I}_3') \le L^{-1}.$$

On the other hand, note that

$$
\begin{aligned}
T_2 = P\Bigg\{ &\sum_{(j,i)\in\mathcal{I}_4'} [(S_{j,i}^2 - \xi_{j,i}^2 - Ln^{-1}) + b_3 Ln^{-1} I(S_{j,i}^2 > \lambda_* Ln^{-1})] \\
&< -(b_1 + b_2)N^{1/2}n^{-1} \\
&+ \sum_{(j,i)\in\mathcal{I}_4'} [(S_{j,i}^2 + n^{-1}\chi_{j,i}^2 - \xi_{j,i}^2 - Ln^{-1})I(S_{j,i}^2 > \lambda_* Ln^{-1})] \Bigg\} \\
\leq P\Bigg\{ &\sum_{(j,i)\in\mathcal{I}_4'} (S_{j,i}^2 - \xi_{j,i}^2 - Ln^{-1}) < -(b_1 + b_2)N^{1/2}n^{-1} \Bigg\} \\
+ P\Bigg\{ &\sum_{(j,i)\in\mathcal{I}_4'} (S_{j,i}^2 + n^{-1}\chi_{j,i}^2 - \xi_{j,i}^2 - Ln^{-1})I(S_{j,i}^2 > \lambda_* Ln^{-1}) \\
&> \sum_{(j,i)\in\mathcal{I}_4'} b_3 Ln^{-1} I(S_{j,i}^2 > \lambda_* Ln^{-1}) \Bigg\} \\
\equiv\ &T_{21} + T_{22}.
\end{aligned}
$$

Set $\tilde{Z}_{j,i} = \xi_{j,i}^{-1} \sum_{(j,k)\in B_i^j} \theta_{j,k} z_{j,k}$. Then $\tilde{Z}_{j,i}$ is a standard Normal random variable and

$$
\begin{aligned}
T_{21} = P\Bigg\{ &\sum_{(j,i)\in\mathcal{I}_4'} (2n^{-1/2}\xi_{j,i}\tilde{Z}_{j,i} + n^{-1}\chi_{j,i}^2 - Ln^{-1}) < -(b_1 + b_2)N^{1/2}n^{-1} \Bigg\} \\
\leq P\Bigg\{ &\sum_{(j,i)\in\mathcal{I}_4'} \chi_{j,i}^2 < -b_1 N^{1/2} + \mathrm{Card}(\mathcal{I}_4')L \Bigg\} \\
+ P\Bigg\{ &\sum_{(j,i)\in\mathcal{I}_4'} \xi_{j,i}\tilde{Z}_{j,i} < -\tfrac{1}{2}b_2 N^{1/2}n^{-1/2} \Bigg\}.
\end{aligned}
$$

If $\mathrm{Card}(\mathcal{I}_4')L \leq b_1 N^{1/2}$, then $P\{\sum_{(j,i)\in\mathcal{I}_4'} \chi_{j,i}^2 < -b_1 N^{1/2} + \mathrm{Card}(\mathcal{I}_4')L\} = 0$. When $\mathrm{Card}(\mathcal{I}_4')L > b_1 N^{1/2}$, equation (45) with $m = \mathrm{Card}(\mathcal{I}_4')L \leq N$ and $d = b_1 N^{1/2}/m$ yields that

$$
P\Bigg\{ \sum_{(j,i)\in\mathcal{I}_4'} \chi_{j,i}^2 < -b_1 N^{1/2} + \mathrm{Card}(\mathcal{I}_4')L \Bigg\} \leq e^{(-1/4)d^2 m} \leq e^{(-1/4)b_1^2} = \frac{\alpha}{2}.
$$

On the other hand, note that $\sum_{(j,i)\in\mathcal{I}_4'} \xi_{j,i}^2 \leq NL^{-1} \cdot 4\lambda_* Ln^{-1} = 4\lambda_* N n^{-1}$ and hence $P\{\sum_{(j,i)\in\mathcal{I}_4'} \xi_{j,i}\tilde{Z}_{j,i} < -\tfrac{1}{2}b_2 N^{1/2}n^{-1/2}\} \leq P(Z < -\tfrac{1}{4}b_2\lambda_*^{-1/2}) \leq \frac{\alpha}{2}$ where $Z \sim N(0,1)$.

We now turn to the term $T_{22}$. Note that $\xi_{j,i}^2 \le 4\lambda_* L n^{-1}$ for $(j,i) \in \mathcal{I}_4'$. Hence

$$T_{22} \le \sum_{(j,i) \in \mathcal{I}_4'} P(2n^{-1/2}\xi_{j,i}\tilde{Z}_{j,i} + 2n^{-1}\chi_{j,i}^2 > (b_3+1)Ln^{-1})$$

$$\le \sum_{(j,i) \in \mathcal{I}_4'} P(\tilde{Z}_{j,i} > \tfrac{1}{2}\xi_{j,i}^{-1}((b_3 - 2\lambda_* + 1)Ln^{-1/2})) + \sum_{(j,i) \in \mathcal{I}_4'} P(\chi_{j,i}^2 > \lambda_* L)$$

$$\le \sum_{(j,i) \in \mathcal{I}_4'} P(\tilde{Z}_{j,i} > 2(\log n)^{1/2}) + \sum_{(j,i) \in \mathcal{I}_4'} \tfrac{1}{2}n^{-2} \le LNn^{-2} \le L^{-1}.$$

Hence, $C(\theta) \le T_1 + T_{21} + T_{22} \le \alpha + 2L^{-1} = \alpha + 2(\log n)^{-1}$.

## REFERENCES

[1] BARAUD, Y. (2004). Confidence balls in Gaussian regression. *Ann. Statist.* **32** 528–551. MR2060168

[2] BERAN, R. (1996). Confidence sets centered at $C_p$-estimators. *Ann. Inst. Statist. Math.* **48** 1–15. MR1392512

[3] BERAN, R. and DÜMBGEN, L. (1998). Modulation of estimators and confidence sets. *Ann. Statist.* **26** 1826–1856. MR1673280

[4] CAI, T. (1999). Adaptive wavelet estimation: A block thresholding and oracle inequality approach. *Ann. Statist.* **27** 898–924. MR1724035

[5] CAI, T. (2002). On block thresholding in wavelet regression: Adaptivity, block size, and threshold level. *Statist. Sinica* **12** 1241–1273. MR1947074

[6] CAI, T. and LOW, M. (2004). Adaptive confidence balls. Technical report, Dept. Statistics, Univ. Pennsylvania.

[7] DONOHO, D. L. and JOHNSTONE, I. M. (1998). Minimax estimation via wavelet shrinkage. *Ann. Statist.* **26** 879–921. MR1635414

[8] GENOVESE, C. R. and WASSERMAN, L. (2005). Confidence sets for nonparametric wavelet regression. *Ann. Statist.* **33** 698–729. MR2163157

[9] HOFFMANN, M. and LEPSKI, O. (2002). Random rates in anisotropic regression (with discussion). *Ann. Statist.* **30** 325–396. MR1902892

[10] JUDITSKY, A. and LAMBERT-LACROIX, S. (2003). Nonparametric confidence set estimation. *Math. Methods Statist.* **12** 410–428. MR2054156

[11] LI, K.-C. (1989). Honest confidence regions for nonparametric regression. *Ann. Statist.* **17** 1001–1008. MR1015135

[12] ROBINS, J. and VAN DER VAART, A. (2006). Adaptive nonparametric confidence sets. *Ann. Statist.* **34** 229–253.

DEPARTMENT OF STATISTICS
THE WHARTON SCHOOL
UNIVERSITY OF PENNSYLVANIA
PHILADELPHIA, PENNSYLVANIA 19104-6340
USA
E-MAIL: tcai@wharton.upenn.edu
        lowm@wharton.upenn.edu